\documentclass[preprint,a4paper,preprint,12pt]{elsarticle}
\usepackage[top=3cm,right=2cm,bottom=2cm,left=3cm]{geometry}
\usepackage{amssymb}
\usepackage{amsmath}
\usepackage{amsthm}
\usepackage{bm}
\usepackage{enumerate}
\usepackage{color}
\usepackage{multirow}
\usepackage{booktabs}
\usepackage[table,xcdraw]{xcolor}
\usepackage{ulem}
\usepackage{graphicx}
\usepackage{lineno}
\usepackage{tabu,hhline,multirow}
\usepackage[ruled,vlined]{algorithm2e}

\begin{document}
	
	\begin{frontmatter}
		
\title{A methodology to evaluate corroded RC structures using a probabilistic damage approach}

\author[fec]{Karolinne O. Coelho \corref{cor1}}

\author[USP]{Edson D. Leonel}
\ead{edleonel@sc.usp.br}

\author[UNILA]{Julio Fl\'orez-L\'opez}
\ead{julio.lopez@unila.br}

\address[fec]{ FEC - Universidade Estadual de  Campinas, R. Josiah Willard Gibbs 85 - Cidade Universitária, Campinas, SP, CEP 13083-839, Brazil}

\address[USP]{Department of Structural Engineering, S\~ao Carlos School of Engineering, University of S\~ao Paulo, Av. Trabalhador S\~ao Carlense, 400, S\~ao Carlos, 13566-590, Brazil}

\address[UNILA]{ Latin-American Institute of Technology, Infrastructure and Territory, Federal University of Latin-American Integration, Av. Silvio Américo Sasdelli, 1842, Foz do Igua\,cu, 85866-000, Brazil}

\cortext[cor1]{Corresponding author - k192247@dac.unicamp.br}

\begin{abstract}
Several aspects influence corrosive processes in RC structures, such as environmental conditions, structural geometry, and mechanical properties. Since these aspects present large randomnesses, probabilistic models allow a more accurate description of the corrosive phenomena. On the other hand, the definition of limit states, applied in the reliability assessment, requires a proper mechanical model. In this context, this study proposes an accurate methodology for the mechanical-probabilistic modelling of RC structures subjected to reinforcements’ corrosion. To this purpose, an improved damage approach is proposed to define the limit states for the probabilistic modelling, considering three main degradation phenomena: concrete cracking, rebar yielding, and rebar corrosion caused either by chlorides or carbonation process. The stochastic analysis is evaluated by the Monte Carlo simulation method due to the computational efficiency of the LDMC. The proposed mechanical-probabilistic methodology is implemented in a computational framework and applied to the analysis of a simply supported RC beam, and a 2D RC frame. Curves illustrate the probability of failure over a service life of 50 years. Moreover, the proposed model allows drawing the probability of failure map and then identify the critical failure path for progressive collapse analysis. Collapse path changes caused by the corrosion phenomena are observed.
\end{abstract}

\end{frontmatter}

\section{Introduction}
The modelling of corrosion in RC structures is usually divided into two stages\citep{tuutti1982corrosion}: the initiation and the propagation phases. The initiation period consists of modelling the diffusion agents ($\mbox{CO}_2$ or $\mbox{Cl}^-$) mechanisms into the concrete pores, whist the propagation is associated with the reinforcements' mass loss and mechanical degradation properties over time.

The initiation phase is a complex topic in which the definition of the variables that influence this phase is not a consensus. Since the steel-concrete interface (SCI) is one of the main reason of starting corrosive processes, a recent study by \citet{angst2019critical} focused on the variables that influence the characteristics of the SCI for chloride-induced corrosion. The authors verified that several variables affect the SCI, such as reinforcing steel type, steel surface preparation, and factors related to the concrete micro and macrostructure, for instance, $w/b$ ratio, cement type, and pre-existing cracks. Moreover, environmental factors lead to the reduction of the initiation time, such as moisture and temperature. These aspects influence not only the SCI but also the concrete permeability and, consequently, the diffusion of the aggressive agents into the concrete \citep{shafikhani2019quantification,shah2018carbonation}.

For simplicity, this study applies analytical models for diffusion. For carbonation, the modelling follows a semi-empirical approach based on the mass balance equation for $\mbox{CO}_2$ reactions \citet{papadakis1992effect}. On the other hand, for the chloride case, the analytical modelling is based on Fick's law for diffusion that simulates the flux of $\mbox{Cl}^-$ ions into a porous media along time. For this formulation, the diffusion coefficient is computed according to \citet{papadakis1996mathematical}. It is worth mentioning that these approaches for diffusion incorporate inherent simplifications associated with the flux direction and saturated conditions into the media. Other formulations can be performed to consider complex and general diffusion problems and using numerical methods such as the Finite Element method \citep{xiao2012fem,benkemoun2017embedded} and the Boundary Element method \citep{pellizzer2020time, chen2017coupled}.

Once defined the initiation time, the propagation stage begins with the destruction of the passivating layer. Then, the corrosive reactions accelerate, leading to reinforcement's mass loss, the rebar's adherence loss, cover cracking and even concrete spalling. To access how fast the structure is deteriorating, the corrosion rate is the variable that indicates the velocity of the corrosive reaction. This variable depends on several parameters such as concrete cover, water/cement ($w/c$) ratio, temperature, and moisture \citep{angst2009critical}. Some approaches are described in the literature for modelling this phase, as presented by \citet{vennesland2013recommendation}, \citet{kiani2012response}, \citet{xu2018numerical}, \citet{yu2021integrated}.

As a result of the propagation of corrosive reactions, the mechanical capacity and stiffness decrease over time. Given the complexity of the phenomena, substantial efforts have been made toward accurate and robust approaches for handling the mechanical effects caused by corrosion in RC structures. Due to the large randomness of such a phenomenon, several of these studies apply stochastic models to predict the probability of failure of corroded RC structures \citep{estes1998relsys, bastidas2018reliability,huang2020stochastic, jia2020stochastic}. However, many of these models use deterministic approaches based on analytical/experimental data. Others use experimental data as the limit state equation for probabilistic analysis. In spite of the mentioned advances, the mechanical modelling of the entire structural life combining accurate numerical methods and stochastic analysis are scarce in the literature, which justifies the development of the present research \citep{liberati2014nonlinear}.

In this regard, this study presents a robust numerical methodology for the mechanical-probabilistic modelling of RC structures subjected to reinforcements’ corrosion with low computational costs, besides the complexity of the corrosive processes in RC structures. The proposed methodology, based on the lumped damage mechanics (LDM) theory, includes the corrosion evolution law to the constitutive model, along with the damage and plasticity. Furthermore, this framework is coupled with a stochastic model to contemplate reliability analysis and perform a full structural simulation for a service life of 50 years. This methodology can be easily implemented in computational software and is flexible regarding the corrosive laws for diffusion and propagation.

The LDM theory, the mechanical model applied in this study, is a well-known model that couples the concepts of damage and fracture mechanics with the plastic hinge idea \citep{florez1998frame, santoro2013damage}. Thus, the plastic hinges incorporate the damage index, which is named generally as the inelastic hinge. Consequently, the LDM approach concentrates the inelastic effects at the element ends. This simple idea provides an extensive reduction in computational time consuming, which makes this approach computationally efficient. Several studies demonstrated the accuracy of the LDM in complex structural analysis, such as 3D frames, concrete arches, high cycle fatigue and cyclic or impact loadings \citep{marante2003three, amorim2014simplified, bai2016macromodeling, barrios2020numerical}.

The LDMC consider that the mechanical effects associated with the corrosion are incorporated into a new state variable named as the corrosion state variable. Such a variable is defined as the ratio between the corroded and the non-corroded rebar cross-section area. The corroded area is calculated according to the corrosion type, pitting or uniform, as presented by \cite{stewart2004spatial}. Thus, the LDMC includes a corrosion evolution law based on the experimental observations. In this study, this law is given by the corrosion rate equations presented by \cite{goctermann2000dura}, for the carbonation case, and by \cite{vu2000structural}, for the chloride case. It is worth mentioning that the corrosive laws applied in this study are only a guide. Any other formulation can be coupled to the LDMC, which makes this model very generic and robust.

Another advantage of this mechanical modelling is that it provides the damage value for each hinge and the collapse time. Then the limit state equation is written as a function of a threshold damage value for the stochastic analysis. Due to the low computational cost of the LDMC, the simulation is performed by the direct coupling of LDMC and the Monte Carlo simulation algorithm. Consequently, the probability of failure curves along time, for both global and local perspectives, are achieved with good precision. Stochastic processes based on the extreme value theory represent the external loading, which makes realistic the loading modelling. Since the proposed formulation provides the individual probability of failure per hinge, the critical failure path is faster assessed.

Two applications are presented in this study: a simply supported RC beam and a 2D hyperstatic RC frame. In both cases, the global and the individual (per hinge) probabilities of failure over time were obtained successfully. In particular, for the hyperstatic application, the corrosion processes change the critical collapse path over time. Therefore the critical collapse path predicted during the design phase, i.e. without corrosion scenario, is considerably different from the collapse configuration observed when the corrosion is accounted. As a result, a failure mode different from the predicted without the corrosion scenario occurs.

\section{Mechanical model for corrosion}
\subsection{The Lumped Damage Model for Corrosion (LDMC)}

The lumped damage theory is based on adding the fracture mechanics and classical damage concepts into plastic hinges. In such a theory, the damage variable is coupled to the plastic hinges, which are generally named as inelastic hinges. The concentration of the damage index at the hinges provides large computational efficiency. Despite being simpler, this approach is efficient and provides results as accurate as those observed by complex and refined damage models.

The present study extends the LDM to contemplate the corrosive mechanisms. This improved formulation leads to the Lumped Damage Model for Corrosion (LDMC). Therefore, this study proposes a numerical formulation capable of accomplishing the inelastic mechanical modelling of RC structures subjected to reinforcements' corrosion. The main degradation phenomena are accounted: concrete cracking, rebar yielding, and cross-section reduction and yield stress reduction caused by corrosion. Such improvements represent the contributions of the present study.

\begin{figure}[h]
	\centering
	\includegraphics[scale=0.9]{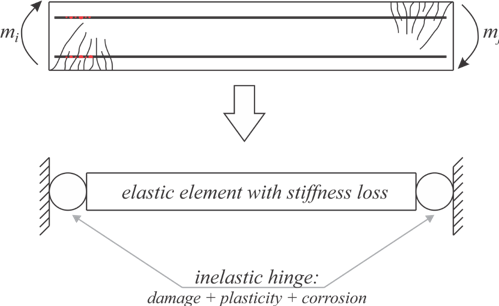}
	\caption{A RC element and its respective LDMC model.}
	\label{fig:rcelement}
\end{figure}

The Fig. \ref{fig:rcelement} illustrates an RC element subjected to bending moments $m_i$ and $m_j$ applied at its ends $i$ and $j$, respectively. The LDMC assumes the mechanical degradation phenomena concentrated at the element ends. The plastic hinge represents the reinforcements' yielding, whereas the concrete cracking is described by the damage variable added to the inelastic hinge. Furthermore, the corrosion state variable is added to the inelastic hinges.

As previously mentioned, the proposed formulation concentrates the corrosion effects at the inelastic hinges. This is a proper hypothesis because the concrete cracking increases the porosity at this region. Consequently, the coefficient of diffusion at the inelastic hinge increases as well as the corrosion rate. Moreover, the formulation assumes that the corrosive processes reduce the structural stiffness. As a result, both flexural and the axial element stiffness are penalized.

\begin{figure}[h]
	\centering
	\includegraphics[scale=0.9]{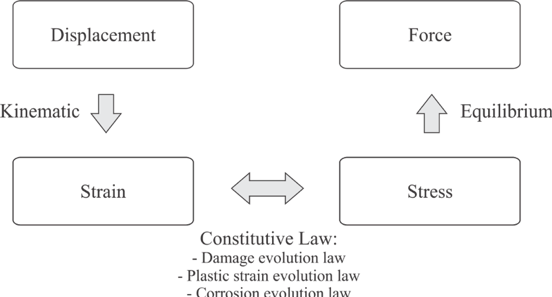}
	\caption{The flowchart for the LDMC.}
	\label{fig:flowchart}
\end{figure}

Kinematic, equilibrium and constitutive laws enable the assessment of stresses, strains, damage index and the relative rotation on the structural elements, as presented by the flowchart illustrated in Fig. \ref{fig:flowchart}.

The constitutive law incorporates three other laws related to the inelastic effects represented by the formulation: the damage evolution law, the plastic strain evolution law and the corrosion evolution law. The equations utilized for representing each of these effects are described in the following.

\subsubsection{Kinematic and Equilibrium law}

The mechanical responses are achieved with the FEM by the planar frame element. This element supports six degrees of freedom (three for each node). The generalized displacement vector is given by Eq. \ref{eq:u}.

\begin{equation}
	\mathbf{u} = \begin{Bmatrix}u_i & w_i & \theta_i & u_j & w_j & \theta_j \end{Bmatrix}^T
	\label{eq:u},
\end{equation}

\noindent in which $u$, $w$ and $\theta$ are the degrees of freedom: the horizontal and vertical displacements and the rotation, respectively.

The displacements of the finite element are related to its strains by the kinematic equations, as presented the Eq. \ref{eq:kinematics}.

\begin{equation}
	\Delta\boldsymbol{\phi} = \mathbf{B}_0\Delta\mathbf{u}
	\label{eq:kinematics}
\end{equation}

As well known, the variation on the generalized displacement vector is associated with the variation on the strain ($\Delta\Phi$). This relationship is governed by the kinematic transformation matrix, presented in Eq. \ref{eq:kinematicmatrix} for the planar frame element.

\begin{equation}
	\mathbf{B}_0 =
	\begin{bmatrix}
		\frac{\sin{(\alpha)}}{L} & -\frac{\cos{(\alpha)}}{L} & 1 & -\frac{\sin{(\alpha)}}{L} & \frac{\cos{(\alpha)}}{L} & 0\\
		\frac{\sin{(\alpha)}}{L} & -\frac{\cos{(\alpha)}}{L} & 0 & -\frac{\sin{(\alpha)}}{L} & \frac{\cos{(\alpha)}}{L} & 1\\
		-\cos{(\alpha)} & -\sin{(\alpha)} & 0 & \cos{(\alpha)} & \sin{(\alpha)} & 0
	\end{bmatrix}
	\label{eq:kinematicmatrix}
\end{equation}

\noindent where $\alpha$ is the angle between the finite element and the coordinate reference system and $L$ is the length of the finite element.

For the planar frame element, the generalized strain vector is composed of the relative rotations at the element ends ($\phi_i$ and $\phi_j$, respectively) and the total elongation, as presented in Eq. \ref{eq:elongation}.

\begin{equation}
	\boldsymbol{\phi} = \begin{Bmatrix}\phi_i & \phi_j & \delta \end{Bmatrix}^T
	\label{eq:elongation}.
\end{equation}

The well-known equilibrium law associates the external forces vector $\mathbf{p}$ to the generalized stress vector $\mathbf{m}$ as follows

\begin{equation}
	\mathbf{B}_0^T\mathbf{m} = \mathbf{p}
\end{equation}

In the present study, only the quasi-static loading condition is considered. In addition, the hypothesis of small displacements and strains are assumed.

\subsection{Constitutive law}

The constitutive law relates the generalized strain to the stress vector, as presented by Eq. \ref{eq:constitutivelaw}.

\begin{equation}
	\boldsymbol{\phi} - \boldsymbol{\phi}_p = \mathbf{F}(\mathbf{d},c)\mathbf{m} + \boldsymbol{\phi}_0
	\label{eq:constitutivelaw}
\end{equation}

\noindent where $\boldsymbol{\phi}_0$ is the initial strain vector (nil for isostatic structures) and $\mathbf{F}(\mathbf{d},c)$ is the flexibility matrix, which depends on the damage values $\mathbf{d}$ and the corrosion state variable $c$.

The hypothesis of strain equivalence enables writing the flexibility matrix for a given element as follows:

\begin{equation}
	\mathbf{F}(\mathbf{d},c) = \mathbf{F}_0 + \mathbf{C}(\mathbf{d},c)
	\label{eq:flexibilitymatrix},
\end{equation}

\noindent in which $\mathbf{F}_0$ is the flexibility matrix for the elastic element and $\mathbf{C}(\mathbf{d},c)$ is the additional flexibility caused by the cracking and the corrosive processes. Thus, the flexibility matrix of a damaged element is evaluated as follows:

\begin{equation}
	\mathbf{F}(\mathbf{d},c) = \begin{pmatrix}
		\frac{L}{3EI(c)(1-d_i)} & -\frac{L}{6EI(c)} & 0\\
		-\frac{L}{6EI(c)} & \frac{L}{3EI(c)(1-d_j)} & 0\\
		0 & 0 & \frac{L}{AE(c)}
	\end{pmatrix}
	\label{eq:flexibilitymatrix},
\end{equation}

\noindent where $EI(c)$ and $AE(c)$ are the flexural and the axial stiffness, respectively, function of the corrosion variable $c$.

As presented in Eq. \ref{eq:flexibilitymatrix}, the damage index penalizes solely the flexural element stiffness, whereas the corrosion state variable penalizes both flexural and axial element stiffness. Moreover, for time instants smaller than the time for corrosion initiation, the corrosion state variable assumes nil value. Thus, in such case, the flexibility matrix coincides with the elastic condition. The inelastic phenomena are represented by penalizing the constitutive law and, consequently, the element stiffness. Thus, the damage law, the plastic strain evolution law and the corrosion evolution law complement the constitutive equations.

\subsubsection{The corrosion evolution law}

The corrosion evolution modelling is accomplished by adding the corrosion state variable to the inelastic hinge. This new variable is defined as follows:

\begin{equation}
	c = \frac{A_c}{A_0}
	\label{eq:corrosionvariable}
\end{equation}

\noindent where $A_c$ is the corroded area and $A_0$ indicates the initial rebar cross-section area.

Therefore, the corrosion state variable varies between 0, in the case of non-corroded reinforcement, and 1, for the entire reinforcement's deterioration. The cross-section area is defined in terms of the corrosion state variable as follows:

\begin{equation}
	A_s = (1-c)A_0
	\label{eq:steelarea}.
\end{equation}

The corrosion evolution law is the variation of the rebar cross-section area concerning the diameter penalization. For the chloride corrosion, the reinforcement’s cross-section area varies as a function of the pit depth, $p$ (pit corrosion case). For carbonation corrosion, the rebar’s cross-section area is reduced as a function of the diameter ($\varphi$) variation (uniform corrosion case). Both equations are described in Eq. \ref{eq:ccarb} and \ref{eq:cchl} for carbonation and chloride corrosion, respectively.

\begin{equation}
	c_{carb} = \frac{1}{A_0}\frac{\mbox{d}A_c}{\mbox{d}\Delta\varphi}\Delta\varphi(i_{corr},x_c)
	\label{eq:ccarb}
\end{equation}

\begin{equation}
	c_{chl} = \frac{1}{A_0}\frac{\mbox{d}A_c}{\mbox{d}p}p(i_{corr},C_0,C_{lim})
	\label{eq:cchl}
\end{equation}

As previously mentioned, the carbonation leads to uniform corrosion. Consequently, the modelling of this phenomenon often assumes the uniform reduction on the rebar’s diameter along the surface reached by the carbonation front. In the proposed formulation, the diameter evolution proposed by \cite{vennesland2013recommendation} was utilized:

\begin{equation}
	\Delta\varphi = 0.0232i_{corr}(t-t_{ini})
	\label{eq:deltaphi}.
\end{equation}

On the other hand, the chloride corrosion is associated with pitting corrosion. The pit depth is a function of the corrosion time initiation, as observed in Eq. \ref{eq:pit}.

\begin{equation}
	p = 0.0116i_{corr}R(t-t_{ini})
	\label{eq:pit},
\end{equation}

\noindent where $R$ indicates the relationship between the maximum and mean pit depth (with $p$ in millimetres). As mentioned by \cite{vu2000structural}, this parameter is probabilistic, with mean value equal to 5.08. Because of the complexity involved in determining its statistical properties (it depends on corrosion time initiation, chloride content, porosity, among others), in this study, $R$ was assumed as 5.08.

The rebar’s cross section area reduces as a function of the pit depth. This reduction is calculated from the Eq. \ref{eq:A}, \ref{eq:A1}, \ref{eq:A2}, as proposed by \cite{val1998effect}.

\begin{equation}
	A = \begin{cases}
		\frac{\pi\varphi^2}{4}-A_1-A_2, & p\leq\frac{\sqrt{2}}{2}\varphi\\
		A_1-A_2, & \frac{\sqrt{2}}{2}\varphi<p<\varphi\\
		0, & p>\varphi
	\end{cases}
	\label{eq:A}
\end{equation}

\begin{equation}
	A_1 = \frac{1}{2}\left[\theta_1\left(\frac{\varphi}{2}\right)^2 - a\left|\frac{\varphi}{2}-\frac{p^2}{\varphi}\right|\right]
	\label{eq:A1}
\end{equation}

\begin{equation}
	A_2 = \frac{1}{2}\left[\theta_2p^2-a\frac{p^2}{\varphi}\right]
	\label{eq:A2}
\end{equation}
where $a, \theta_1, \theta_2$ are geometric parameters assessed as follows:

\begin{equation}
	a = 2p\sqrt{q-\left(\frac{p}{\varphi}\right)^2}
	\label{eq:a}
\end{equation}

\begin{equation}
	\theta_1 = 2\arcsin{\left(\frac{2a}{\varphi}\right)}
	\label{eq:theta1}
\end{equation}

\begin{equation}
	\theta_2 = 2 \arcsin{\left(\frac{a}{\varphi}\right)}
	\label{eq:theta2}
\end{equation}

Despite the formulation is based on empirical/analytical approaches, other equations may be utilized for assessing/improving the corrosion evolution law. Therefore, any other corrosion model can be coupled with the proposed formulation, which illustrates its flexible and easily updated aspects. 

It is worth mentioning that the corrosion evolution law modifies solely the constitutive relationships. Consequently, it affects the nonlinear mechanical-material behaviour by penalizing the element stiffness and by enlarging the damage variable. Thus, the latter is composed of material and corrosion components.

\subsubsection{The damage evolution law}

The damage evolution law is based on the energy balance criterion formulated by \citet{griffith1921vi}. Griffith affirms that cracks only propagate when the available energy reaches the energy value required to the crack extension. Thus, the total energy reduces/stabilize because of the energy dissipation, which leads to the additional crack surfaces.

Therefore, the damage evolution law is proposed through the comparison of the energy release rate and the cracking resistance at the inelastic hinge. The relationship defined by Griffith is mathematically expressed by the Eq. \ref{eq:damagelaw1} and \ref{eq:damagelaw2}.

\begin{equation}
	\begin{cases}
		\Delta d_i = 0, & \mbox{if }G_i<R_i\\
		G_i=R_i, & \mbox{if }\Delta d_i>0
	\end{cases}
	\label{eq:damagelaw1}
\end{equation}

\begin{equation}
	\begin{cases}
		\Delta d_j = 0, & \mbox{if }G_j<R_j\\
		G_j=R_j, & \mbox{if }\Delta d_j>0
	\end{cases}
	\label{eq:damagelaw2}
\end{equation}

\noindent where the sub-indexes $i$ and $j$ are related to the inelastic hinges $i$ and $j$, respectively.

The energy release rate is based on the derivative of the complementary energy with respect to the damage parameters. When the frame element with inelastic hinges is accounted, the complementary energy is defined as the generalized stress work with respect to the generalized strain \cite{florez1998frame}. Therefore, the damage evolution is defined in Eq. \ref{eq:gi} and \ref{eq:gj} for the hinges $i$ and $j$, respectively, as follows:

\begin{equation}
	G_i = \frac{\partial W_b}{\partial d_i} = \frac{L m_i^2}{6EI(c(1-d_i)}
	\label{eq:gi}
\end{equation}

\begin{equation}
	G_j = \frac{\partial W_b}{\partial d_j} = \frac{L m_j^2}{6EI(c)(1-d_j)}
	\label{eq:gj}.
\end{equation}

\begin{figure}[h]
	\centering
	\includegraphics[scale=0.6]{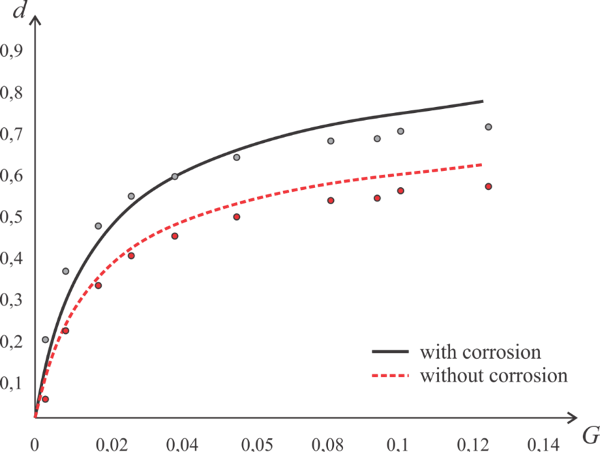}
	\caption{RC cracking resistance curve based on the energy release rate – with and without corrosion.}
	\label{fig:rcresistancecurve}
\end{figure}

The Griffith statement, the basis of the LDM cracking resistance function, is a semi-empirical equation obtained through the stiffness variation method. Additional details of the experimental tests are presented in \cite{florez1998frame}. The corresponding damage and bending moment values, obtained from the experimental tests, are utilized into the energy release rate equation, which is adjusted into a curve, as showed in Fig. \ref{fig:rcresistancecurve}. For RC structures subjected to reinforcements’ corrosion, the energy release rates in corroded condition are smaller than the release rates in non-corroded condition. Thus, during the corroded condition, the concrete cracking process is triggered by a smaller energy amount.

During the crack propagation, the energy release rate is equal to the crack resistance function. Therefore, the curve illustrated in Fig. \ref{fig:rcresistancecurve} is fitted as a function of three variables: $R_0$, $q$, and $d$, as proposed in Eq. \ref{eq:rd}.

\begin{equation}
	R(d) = R_0 + q\frac{\ln{(1-d)}}{1-d}
	\label{eq:rd}
\end{equation}

\begin{figure}
	\centering
	\includegraphics[scale=0.6]{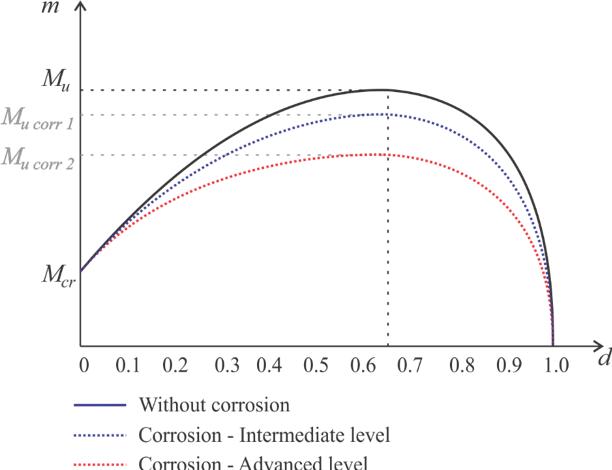}
	\caption{Bending moment as a function of damage – with and without corrosion.}
	\label{fig:momentvsdamage}
\end{figure}

The parameters $R_0$ and $q$ depend on the mechanical characteristics of the element. Both parameters are calculated from the bending moment vs. damage index curve. This curve, presented in Fig. \ref{fig:momentvsdamage}, is obtained by matching the energy release rate equation with the cracking resistance. This condition is achieved when the crack propagation occurs, which is mathematically expressed as follows:

\begin{equation}
	m^2 = \frac{6EI(1-d)^2}{L}R_0 + \frac{6qEI}{L}(1-d)\ln{(1-d)}
	\label{eq:m2}
\end{equation}

The $R_0$ parameter is obtained when nil damage index value is observed and the bending moment is equal to the cracking moment. Therefore, this parameter is associated with the initial cracking resistance, as observed by the following equation:

\begin{equation}
	R_0 = \frac{M_{cr}^2L}{6EI}
	\label{eq:r0}.
\end{equation}

The parameter $q$ is expressed as a function of the ultimate bending moment and its respective damage value. The inflexion point at the curve illustrated in Fig. \ref{fig:momentvsdamage} determines the ultimate damage value. Thus, such a point is determined by deriving the function introduced in Eq. \ref{eq:m2} with respect to the damage variable and equalling it to zero.

It is worth mentioning that Fig. \ref{fig:momentvsdamage} also illustrates the evolution of bending moment and damage index for the corroded case. As illustrated in this figure, the corrosive processes lead to the reduction of the ultimate moment. Therefore, threshold values of damage index are observed with smaller values of bending moment in the corroded case. Moreover, the corrosion processes increase the concrete cracking at the hinges because of the reduction on the element resistance.

\subsubsection{The plastic strain evolution law}

The plastic strain evolution law establishes that plastic rotation does not vary if the yield function $f$ is lesser than zero, as observed in Eq. \ref{eq:plasticitylaw1} and \ref{eq:plasticitylaw2}, for the hinges $i$ and $j$, respectively. Thus, plastic strains do not evolve into the elastic domain.

\begin{equation}
	\begin{cases}
		\mbox{d}\phi_{pi} = 0, & \mbox{if } f_i < 0\\
		f_i = 0, & \mbox{if }\mbox{d}\phi_{pi} \neq 0
	\end{cases}
	\label{eq:plasticitylaw1}
\end{equation}

\begin{equation}
	\begin{cases}
		\mbox{d}\phi_{pj} = 0, & \mbox{if } f_j < 0\\
		f_j = 0, & \mbox{if }\mbox{d}\phi_{pj} \neq 0
	\end{cases}
	\label{eq:plasticitylaw2}.
\end{equation}

The yield function $f$ utilized by the proposed formulation includes the effects of damage propagation and kinematic hardening. Thus, $f$ is defined as follows:

\begin{equation}
	f = \left|\frac{m}{1-d} - c_{plast}\phi_p\right| - k_0
\end{equation}

\noindent where $c_{plast}$ and $k_0$ are time and element-dependent constants. The $c_{plast}$ and $k_0$ values are assessed as a function of the yield damage $d_p$ and the plastic moment, through the relationship between bending moment and damage index, previously mentioned in Eq. \ref{eq:m2}. When the yield damage value is reached, the plastic rotation is nil. Consequently, the yield function is also nil. Therefore, one observes that $k_0$ is the effective yield moment, as proposed in Eq. \ref{eq:k0}.

\begin{equation}
	k_0 = \frac{M_p}{1-d_p}
	\label{eq:k0}
\end{equation}

The yield function is also nil when the ultimate bending moment is reached. Thus, based on Eq. \ref{eq:m2}, coefficient $c_{plast}$ is a function of the ultimate plastic rotation as follows:

\begin{equation}
	c_{plast} = \frac{1}{\phi_{pu}}\left(\frac{M_u}{1-d_u}-\frac{M_p}{1-d_p}\right)
	\label{eq:cplast}
\end{equation}

Both $c_{plast}$ and $k_0$ parameters are affected by the corrosive processes. As observed in Fig. \ref{fig:momentvsdamage}, the stiffness loss leads to the decrease of the plastic moment, which implies into rebar yielding for loading values smaller than the planned during the design phase.

\section{Stochastic model}

\subsection{The extreme value theory}

The design codes specify service life of 50 years for usual RC structures. Since assuming a constant load value during 50 years is not an accurate hypothesis, the probabilistic model must account properly for loading modelling during the period. In this regard, the Extreme Value process is an accurate approach for evaluating the load processes over time during a stochastic analysis. In this case, the Extreme Value process is defined, for a cumulative density function (CDF), as follows \citep{beirlant2006statistics}:

\begin{equation}
	F_{Yn}(y) = \left[F_X(y)\right]^n
	\label{eq:cdf}
\end{equation}

\noindent where $n$ represents the number of sample variables. In the proposed formulation, n represents the number of years.

The extreme distribution presented in Eq. \ref{eq:cdf} tends to an asymptotic extreme distribution, which could be classified as Gumbel, Frechet and Weibull (types I, II e III, respectively), according to the original CDF assumed ($F_X(y)$). Several studies, for instance \cite{astroza2019dual} and \cite{oudah2017design}, demonstrate that accidental loading may be properly represented by the Gumbel distribution for maximum, which CDF is defined as follows:

\begin{equation}
	F_X(x) = \exp{\left[-\exp{\left[-\omega\left(x-u_n\right)\right]}\right]}
\end{equation}

\noindent in which $\omega$ indicates the form parameter and $u_n$ the maximum characteristic value given by Eq. \ref{eq:omega} and Eq. \ref{eq:un}. $\mu$ and $\sigma$ are the mean and standard deviation, respectively.

\begin{equation}
	\omega = \frac{\pi}{\sqrt{6}\omega}
	\label{eq:omega}
\end{equation}

\begin{equation}
	u_n = \mu-\frac{0.577216}{\omega}
	\label{eq:un}.
\end{equation}

Bearing in mind that the maximum characteristic of a given extreme distribution is given by Eq. \ref{eq:fxun}, and the design codes define the maximum characteristic value for 50 years, the value $u_n$ for any other year of service life is obtained by Eq. \ref{eq:unu50}.

\begin{equation}
	F_X(u_n) = P\left[\{ X\leq u_n\}\right] = 1-\frac{1}{n}
	\label{eq:fxun}
\end{equation}

\begin{equation}
	u_n = u_{50} + \ln{\left[-\ln{\left(\frac{50-n}{50}\right)}\right]}
	\label{eq:unu50}
\end{equation}

\subsection{The Monte Carlo simulation for the stochastic corrosion analysis}

The Monte Carlo simulation (MCS) method is a well-known numerical simulation technique applied in uncertainty quantification problems. In a brief overview, the MCS assesses the probability of failure from a sampling of random variables. The sampling is constructed for each random variable, considering the corresponding statistical distribution assigned (the CDF curves, for instance). For stochastic simulations, this sampling is assessed at each time step. The MCS computes the probability of failure by simulating the limit state function for each random variable. Then, the structural failure is observed when the sampling points leads to the failure domain. The probability of failure is given by the ratio of the sampling points at the failure domain and the total amount of simulations. Thus:

\begin{equation}
	P_f = \int_{\Omega} I[x]f_X(x)\ \mbox{d}x = \frac{1}{n_t}\sum_{j=1}^{n_t}I[x_j] = \frac{n_f}{n_t}
	\label{eq:pf}
\end{equation}

\noindent where $n_t$ is the sampling range and $n_f$ the amount of observed failures. $I[x]$ is 1 for failure condition and nil for the safe condition.

It is worth mentioning that crude MCS requires a large amount of limit state functions simulations for assessing small values of probability of failure, i.e., $10^{-3}$ or lesser. Because these small values are currently observed in structural applications, this simulation method normally leads to high computational costs. Then, for complex structural systems in which complex mechanical models are involved, the use of MCS may be prohibitive, because complex mechanical models usually require large computational effort.

Nevertheless, this study presents a mechanical-probabilistic approach based on the direct coupling of LDMC and MCS for the simulation of RC structures subjected to reinforcements' corrosion. The robustness, accuracy and computational efficiency of the lumped damage modelling enable the probabilistic analysis through the direct coupling with the MCS \cite{coelho2017mechanical}. This coupled model, based on the improved version of the lumped damage theory, is one contribution of the present study. In addition to the probability of failure, the proposed formulation enables the determination of the progressive collapse modelling accounting for the uncertainties, which is its main contribution.

\begin{figure}[h]
	\centering
	\includegraphics{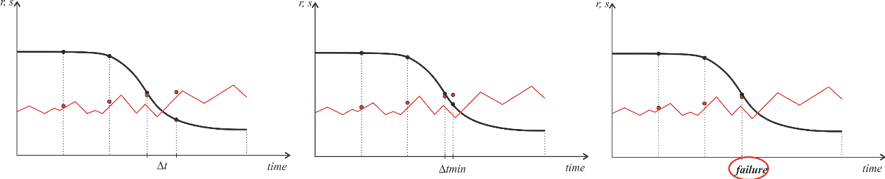}
	\caption{General scheme of the LDMC analysis for a Monte Carlo simulation - in red the solicitation by the Extreme Value theory, and in black the structural capacity (or resistance) decreasing over time due the corrosion and damage phenomena.}
	\label{fig:MCSanalysis}
\end{figure}

The general scheme for the stochastic simulation is presented in Fig. \ref{fig:MCSanalysis}. For each simulation $i$, the mechanical behaviour is modelled into a time step $t$ assuming the resistance decrease caused by the corrosion phenomena and the loading increment calculated from the extreme value theory. If the failure condition is achieved, the time step is reduced to determine the failure time accurately. After $n$ simulations, the probability of failure for each time instant is achieved with the Eq. \ref{eq:pf}.

\section{The algorithm for the mechanical probabilistic modelling of RC corroded structures}

The proposed methodology couples the mechanical and stochastic model to simulate the entire service life of RC structures. The output is the probability of failure, divided into two groups: the probability of failure for the whole structure and for each inelastic hinge (global and local probability of failure, respectively). These data are stored into vectors with 50 position, where each position means the probability of failure in one specific year. An overview of the entire algorithm, implemented using Fortran, is described in Algorithm 1.

The analysis starts by evaluating the initiation time using the diffusion coefficient described in \citet{papadakis1992effect} and \citet{papadakis1996mathematical} for chloride and carbonation corrosion, respectively. Since we want to access the probability of failure due to corrosive processes, the mechanical analysis starts only if the initiation time is less than 50 years. The mechanical analysis computes the global and local collapse time, shown in more details in Algorithm 2.\newline

\begin{algorithm}[H]
	\SetAlgoLined
	\KwResult{Global and local probabilities of failure for each year, over 50 years.}\
	
	read the input file with the probabilistic data and the geometric mesh\;
	initialize the failure vectors with 50 positions fulfilled by 0\;
	\For{$n$ < nsimulations}{
		generate random variables\;
		compute the initiation time for corrosion\;
		\If{initiation time < 50 years}
		{
			run \textbf{Algorithm 2 - LDMC} and get the global and local collapse time for each inelastic hinge\;
			\For{time $\mathrm{from }$ global collapse time $\mathrm{to }$ 50 years}{
				\textit{global failure vector}[\textit{time}]++;
			}
			\ForEach{inelastic hinge}{
				\For{time $\mathrm{from }$ local collapse time $\mathrm{to }$ 50 years}{
					\textit{local failure vector}[\textit{time}]++;
				}
			}
		}
	}
	\ForEach{failure vector}{
		divide each component by \textit{nsimulations}.
	}
	
	\caption{Mechanical probabilistic analysis of RC structures over a service life of 50 years.}\label{alg:fullalgorithm}
\end{algorithm}\vspace{0.5cm}

Since the mechanical analysis is inside the scope of a stochastic problem, the simulation starts by computing the load for a given time, using the Extreme Value theory. Next, the corrosion evolution law is applied to compute the corrosion state variable value ($c$) and then penalize the mechanical and geometric properties, such as the rebar cross-section area and the flexural stiffness. These updated parameters are employed in the inelastic analysis. The LDMC analysis per hinge uses the elastic predictor - inelastic corrector algorithm to identify which hinges are active or not, followed by the assembling of the finite elements and the nonlinear system's solution.

The structural failure is given by the loss of stability characterized by the non-convergence of the nonlinear solver. In case of failure, the algorithm gets the time for the collapse of the whole structure. Even though the structure does not collapse, local failures are observed if the limit state function is violated. Such an equation is written as a function of threshold damage values as follows:

\begin{equation}
	g_2 = d_{thre} - d
\end{equation}

\noindent in which the acceptable damage adopted in this study is 0.5, which, according to \cite{amorim2014simplified}, is the limit value for repair.

Defining the local failures per hinge enables the determination of the critical collapse path, i.e. the most probable failure path – and then, the progressive collapse can be analyzed into a realistic form.\newline

\begin{algorithm}[H]
	\SetAlgoLined
	\KwResult{Global and local structural collapse time.}\
	
	get the input variables: finite element mesh, mechanical properties and initiation time\;
	define the time step\;
	
	\For{time $\mathrm{from }$ initiation time $\mathrm{to}$ 50 years}{
		compute the load using the Extreme Value theory\;
		determine the $c$ value using the corrosion evolution law\;
		penalize the mechanical and geometric properties\;
		\ForEach{hinge}{
			define whether the inelastic hinges are active for plasticity and damage using a elastic predictor - inelastic corrector algorithm\;
		}
		assemble the finite elements\;
		solve the nonlinear system of equations through the Newmark method\;
		\ForEach{inelastic hinge}{
			\If{damage $\geq$ 0.5}{
				get the local collapse time\;
			}
		}
		\If{convergence not found}{
			get the global collapse time\;
			\textbf{break}\;
		}
	}
	
	\caption{Lumped damage mechanics for corrosion.}\label{alg:fullalgorithm}
\end{algorithm}\vspace{0.5cm}

\section{Applications}

Two applications demonstrate the robustness and efficiency of the proposed algorithm. The first application shows the mechanical-probabilistic analysis of a simply supported RC beam. The second one simulates a hyperstatic 2D RC frame. The statistical data utilized in both examples is summarized in Table \ref{tab:statdata}.

\begin{table}[h]
	\centering
	\begin{tabular}{llll}
		\hline
		\hline
		\textbf{}        & \textbf{Mean} & \textbf{COV} & \textbf{PDF} \\ \hline
		$w/c$ ratio        & 0.5           & 0.15         & Normal       \\ \hline
		Moisture  (\%)   & 75            & 0.25         & Normal       \\ \hline
		Temperature (ºC) & 20            & 0.25         & Normal       \\ \hline
		$i_{corr}$ ($\mu A$)       & 0.431         & 0.60         & Lognormal    \\ \hline
		$C_{lim}$             & 0.5           & 0.15         & Normal       \\ \hline
		$C_0$               & 75            & 0.25         & Normal       \\ \hline
		$\mbox{CO}_2$ (\%)         & 2             & 0.25         & Normal       \\ \hline
	\end{tabular}
	\caption{Statistical data for the applications}
	\label{tab:statdata}
\end{table}

It is worth mentioning that \cite{otieno2012prediction} establishes the mean and the standard deviation values for $i_{corr-20}$. The concentration of the corrosive agents into the environment ($\mbox{Cl}^-$ and $\mbox{CO}_2$) are based on the environmental measures performed by \cite{albuquerque2005proposta}, which accounts for two north-east Brazilian cities.

The following deterministic parameters values were assumed: aggregate/cement ratio, 5, and aggregate, cement and water densities of 2560 kg/$m^3$, 2500 kg/$m^3$ and 1000 kg/$m^3$, respectively. In both applications, the amount of $10^6$ Monte Carlo simulations per time step were utilized for assessing the probability of failure. Since the mechanical model is very efficient and fast, the computational time consuming was 11 hours for a serial code, using an ordinary desktop (16 GB RAM and Intel i7-2700K). Even fast, given the high number of simulations, this computational time can be reduced by doing code optimizations or paralellizing the loops - since this code is highly parallelizable code (i.e. the instructions can be executed independently).

\subsection{Simply supported RC beam}

\begin{figure}[h]
	\centering
	\includegraphics[scale=0.6]{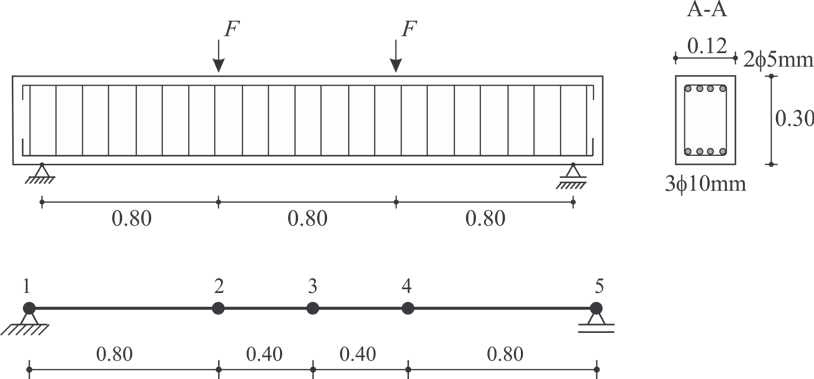}
	\caption{Four-point bending test in an RC beam (dimensions in meters) and the utilized finite element mesh.}
	\label{fig:ex1}
\end{figure}

The first application concerns a simply supported RC beam subjected to two equally distant concentrated loads. The design scheme is illustrated in Fig. \ref{fig:ex1}, as well its respective finite element mesh. As presented in Fig. \ref{fig:ex1}, solely four finite elements are required to accomplish accurately the mechanical-probabilistic analysis. For this mesh, the deterministic lumped damage approach has already demonstrated its accuracy compared to the experimental analysis \citet{coelho2017mechanical}.

The load process is represented by a stochastic approach. The mean load value is assumed as 20 kN, with a coefficient of variation (COV) of 10\%. Table \ref{tab:dataex1} informs the statistical properties for the additional random variables considered in this analysis. It is worth citing that the COV values for concrete cover, compressive resistance and rebar yielding are based on the study presented by \citet{pellizzer2015mechanical}.

\begin{table}[h]
	\centering
	\begin{tabular}{llll}
		\hline
		\hline
		& \textbf{Mean} & \textbf{COV} & \textbf{PDF} \\ \hline
		Concrete cover (cm) & 1.5           & 0.15         & Normal       \\ \hline
		$f_c$ (MPa)            & 38            & 0.10         & Normal       \\ \hline
		$f_y$ (MPa)            & 500           & 0.10         & Lognormal    \\ \hline
		$f_{su}$ (MPa)           & 550           & 0.10         & Lognormal    \\ \hline
		Load (kN)        & 50            & 0.10         & Gumbel Max   \\ \hline
	\end{tabular}
	\caption{Statistical data: simple supported RC beam.}
	\label{tab:dataex1}
\end{table}

\begin{figure}[h]
	\centering
	\includegraphics[scale=0.7]{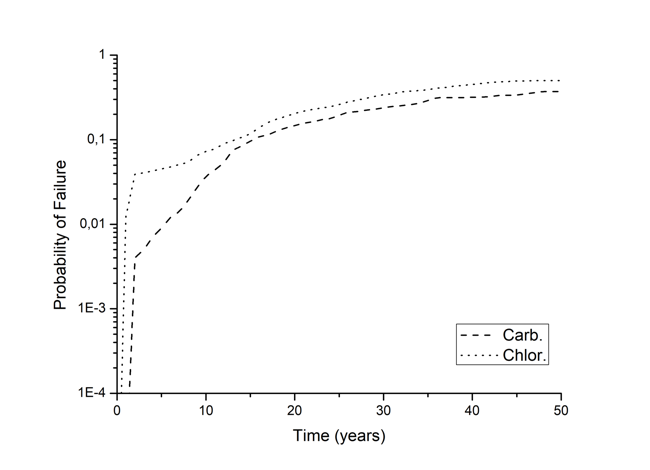}
	\caption{Probability of failure – simply supported RC beam.}
	\label{fig:pfex1}
\end{figure}

Fig. \ref{fig:pfex1} illustrates the evolution of the probability of failure over time. As expected, the corrosion processes increased the probability of structural failure. Due to the corrosive reactions, the low concrete cover and high environmental aggressivity, major values of the probability of failure are observed with only 10 years of structural lifetime. Also, the $w/c$ ratio has an average value smaller than the design codes requirements for the analysed situation, contributing to these high values.

\subsection{2D RC frame}

\begin{figure}[h]
	\centering
	\includegraphics[scale=0.7]{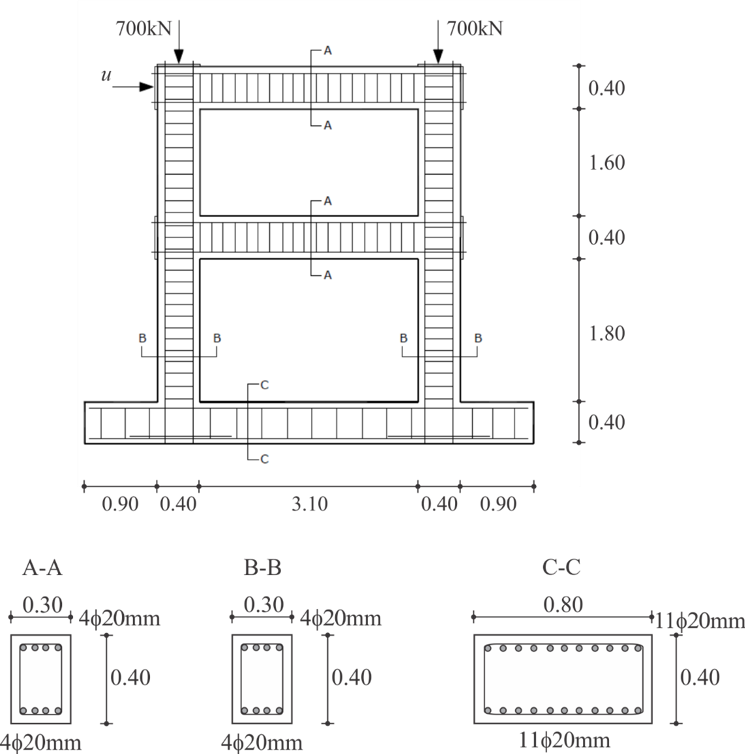}
	\caption{RC 2D frame (dimensions in meters).}
	\label{fig:ex2}
\end{figure}

The second application of this study concerns the hyperstatic 2D RC frame experimentally analysed by \cite{vecchio1986modified} and illustrated in Fig. \ref{fig:ex2}. The loading system is composed of two deterministic vertical loads of 700 kN, which are applied integrally at the analysis begin, and the stochastic horizontal loading for the service life of 50 years. The statistical loading data are presented in Table \ref{tab:dataex2}. This table also presents the statistical data utilized for the additional random variables assumed in this analysis.

\begin{table}[h]
	\centering
	\begin{tabular}{llll}
		\hline
		\hline
		& \textbf{Mean} & \textbf{COV} & \textbf{PDF} \\ \hline
		Concrete cover (cm) & 2.5           & 0.15         & Normal       \\ \hline
		$f_c$ (MPa)            & 30            & 0.10         & Normal       \\ \hline
		$f_y$ (MPa)            & 418           & 0.05         & Lognormal    \\ \hline
		$f_{su}$ (MPa)           & 598           & 0.05         & Lognormal    \\ \hline
		Load (kN)        & 350           & 0.10         & Gumbel Max   \\ \hline
	\end{tabular}
	\caption{Statistical data: 2D RC frame.}
	\label{tab:dataex2}
\end{table}

\begin{figure}[h]
	\centering
	\includegraphics[scale=0.55]{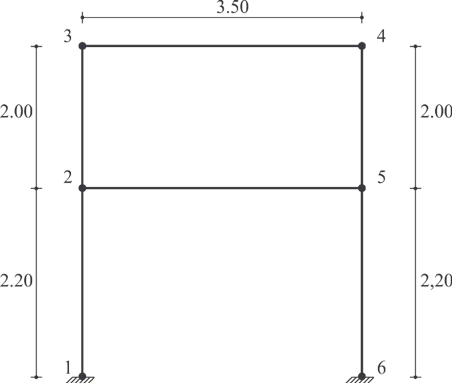}
	\caption{Finite element mesh for the 2D RC frame application (dimensions in meters).}
	\label{fig:meshex2}
\end{figure}

The finite element mesh utilized in mechanical modelling is composed of only six elements, as illustrated in Fig. \ref{fig:meshex2}. This simpler mesh leads to accurate results compared to the experimental curve, as demonstrated by \cite{coelho2017mechanical}. In this hyperstatic structural modelling, the global probability of failure, the local probability of failure (i.e. per hinge) and the critical collapse path are determined for the stochastic loading with mean equal to 200 kN.

\begin{figure}
	\centering
	\includegraphics{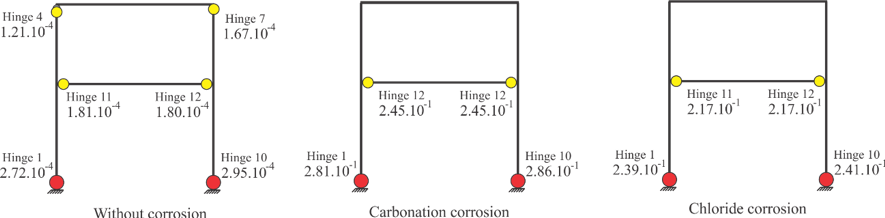}
	\caption{Failure map for the 2D RC frame.}
	\label{fig:failuremap}
\end{figure}

The probability of failure is quantified for the scenarios with or without corrosion, caused either by chloride or carbonation. Fig. \ref{fig:failuremap} illustrates the map of the probability of failure per hinge, with the structural lifetime of 50 years. This figure demonstrates the influence of the corrosion process into probabilistic mechanical behaviour. When the corrosion effects are disregarded, the probability of failure per hinge is of the order of $10^{-4}$. On the other hand, the corrosive processes increased the probability of failure to the order of 10-1 on the hinges 1, 10, 11 and 12. Despite the change into the critical collapse path for chloride and carbonation cases, the maximum local probability of failure is the same for both corrosion types (hinges 1 and 10, at the fixed support).

Since the structural failure occurs by the loss of stability, a proper assessment of the critical collapse path is relevant for probabilistic modelling to evaluate progressive collapse. Moreover, the critical collapse path has significant modifications caused by the corrosion processes, Fig \ref{fig:failuremap}. Without corrosion, the top of the columns (hinges 4 and 7) possess the higher values of the probability of failure. On the other hand, the corrosive reactions lead to failures at the fixed support and the intermediate beam. Consequently, the reinforcements' corrosion may lead to unexpected failure configurations. The maintenance and repair procedures based exclusively in non-corroded cases may be inefficient because the critical hinges are not taken into account.

\begin{figure}[h]
	\centering
	\includegraphics[scale=0.7]{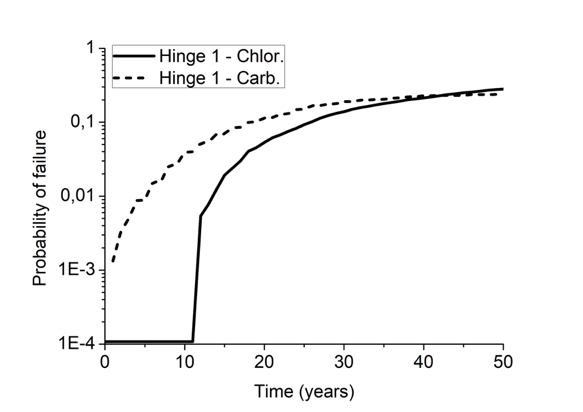}
	\includegraphics[scale=0.7]{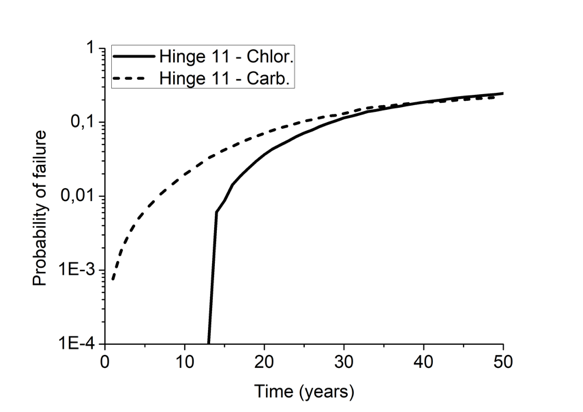}
	\caption{Probability of failure for hinge 1 (left) and 11 (right).}
	\label{fig:pfhinge11}
\end{figure}

\begin{figure}[h]
	\centering
	\includegraphics[scale=0.7]{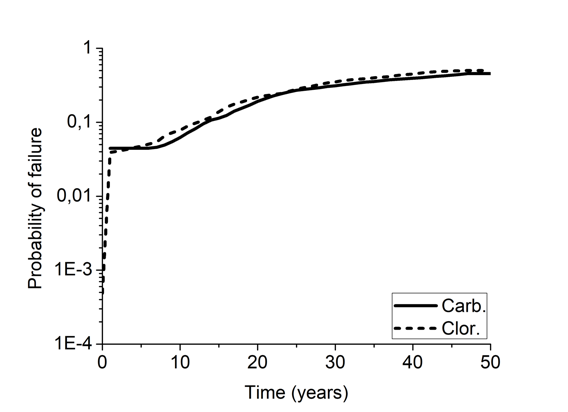}
	\caption{Probability of failure – 2D RC frame.}
	\label{fig:pfex2}
\end{figure}

The probability of failure evolution is illustrated in Fig. \ref{fig:pfhinge11}, for the hinges 1 and 11, and in Figure \ref{fig:pfex2} for the global scope. Since global failure is characterized by the loss of stability, higher values of the probability of failure are observed in the global scope for the corrosion case. In addition, the evolution curves illustrated in Fig. \ref{fig:pfhinge11} and \ref{fig:pfex2} demonstrate different behaviour according to the corrosion type modelled. The chloride corrosion leads to the pitting corrosion scheme, which implies into high corrosion rate localized in a small rebar region. Then, the reinforcement’s mass loss largely increases during the first corrosion years. On the other hand, the carbonation leads to the uniform corrosion scheme, in which the probability of failure growth is smooth along time.\newpage

\section{Conclusions}

This study contributes to the comprehension of the mechanical effects caused by the corrosion in RC structures. The main contribution of this study is the proposition of a methodology to evaluate the structural performance of RC structures subjected to corrosion considering the damage, rebar yielding and stochastic analysis.  For this purpose, the lumped damage model incorporates a new state variable, named as a corrosion state variable, and its respective corrosion evolution law. Such improvement originated the lumped damage model for corrosion (LDMC), which is another contribution of the present study. The LDMC simulates the efforts redistribution during the damage propagation process in RC corroded structures, providing a progressive collapse analysis. As a result, the uncertainty quantification assessment of complex structural systems accounts for the corrosive effects, and thus the critical collapse path is properly evaluated.

The accurate responses were accomplished with the LDMC analysis for isostatic RC structures, taking account the stochastic behaviour for both resistance and loading processes. For hyperstatic structures, the proposed formulation provides the probability of failure map, with an easy visual determination of the critical collapse path. This identification process is especially relevant for the probabilistic-based corrosion analysis because the corrosive processes change the critical collapse path and the local probabilities of failure significantly. The applications demonstrate the potential analysis of the proposed methodology. For instance, it can be applied in studies concerning the progressive collapse of corroded RC buildings. Moreover, high-performance computing schemes, more accurate corrosive models, consideration of concrete spalling, and better statistical data can improve the result's accuracy.

\section*{Acknowledgements}

The financial support from the Brazilian National Council for Scientific and Technological Development (CNPq) is gratefully acknowledged.


\begin{thebibliography}{7}

\bibitem[\protect\citeauthoryear{Albuquerque and Otoch}{2005}]{albuquerque2005proposta}
Albuquerque, A.T. and Otoch, S. (2005), {Proposta de classifica{\c{c}}{\~a}o da agressividade do ambiente na cidade de fortaleza}, \textit{47{\textordmasculine} in Congresso Brasileiro do Concreto}, Recife (in Portuguese).

\bibitem[\protect\citeauthoryear{Amorim \textit{et~al.}}{2014}]{amorim2014simplified}
Amorim, D. L. N. de F., Proenca, S. P. B. and Fl\'orez-L\'opez, J. (2014), {Simplified modeling of cracking in concrete: Application in tunnel linings},  \textit{Engineering Structures}, \textbf{70}, 23--35.

\bibitem[\protect\citeauthoryear{Angst \textit{et~al.}}{2009}]{angst2009critical}
Angst, U., Elsener, B., Larsen, C., and Vennesland, O. (2009), {Critical chloride content in reinforced concrete—a review.}, \textit{Cement And Concrete Research}, \textbf{39}, 1122--1138.

\bibitem[\protect\citeauthoryear{Angst \textit{et~al.}}{2019}]{angst2019critical}
Angst, U.M., Geiker, M.R., Alonso, M.C., Polder, R., Isgor, O.B., Elsener, B., Wong, H., Michel, A., Hornbostel, K., Gehlen, C. and François, R. (2019), {The effect of the steel–concrete interface on chloride-induced corrosion initiation in concrete: a critical review by RILEM TC 262-SCI}, \textit{Materials and Structures}, \textbf{52}(4), 88--113.

\bibitem[\protect\citeauthoryear{Astroza \textit{et~al.}}{2019}]{astroza2019dual}
Astroza, R., Alessandri, Andr{\'e}s and Conte, J. P. (2019), {A dual adaptive filtering approach for nonlinear finite element model updating accounting for modeling uncertainty}, \textit{Mechanical Systems and Signal Processing}, \textbf{115}, 782--800.

\bibitem[\protect\citeauthoryear{Bai \textit{et al.}}{2016}]{bai2016macromodeling}
Bai, Y., Kurata, M., Florez-Lopez, J., and Nakashima, M. (2016), {Macromodeling of crack damage in steel beams subjected to nonstationary low cycle fatigue}, \textit{Journal of Structural Engineering}, \textbf{142}(10),{04016076}.

\bibitem[\protect\citeauthoryear{Barrios and Flórez-López}{2020}]{barrios2020numerical}
Barrios, S. K. M. and Fl{\'o}rez-L{\'o}pez, J. {Numerical quantification of damage reduction in frames retrofitted with FRP bands as bracing elements}, \textit{Engineering Structures}. \textbf{223}, {111178}.

\bibitem[\protect\citeauthoryear{Bastidas-Arteaga}{2018}]{bastidas2018reliability}
Bastidas-Arteaga, E. (2018), {Reliability of reinforced concrete structures subjected to corrosion-fatigue and climate change}, \textit{International journal of concrete structures and materials}, \textbf{12}(1), 10.

\bibitem[\protect\citeauthoryear{Beirlant \textit{et~al.}}{2006}]{beirlant2006statistics}
Beirlant, J., Goegebeur, Y., Segers, J. and Teugels, J. L. (2006), \textit{Statistics of extremes: theory and applications}, John Wiley \& Sons, Chichester, West Sussex, England.

\bibitem[\protect\citeauthoryear{Benkemoun \textit{et al.}}{2017}]{benkemoun2017embedded}
Benkemoun, N, Hammood, M. N., and Amiri, O. (2017), {Embedded finite element formulation for the modeling of chloride diffusion accounting for chloride binding in meso-scale concrete}, \textit{Finite Elements in Analysis and Design}, \textbf{130}, {12--26}.

\bibitem[\protect\citeauthoryear{Chen and  Leung}{2017}]{chen2017coupled}
Chen, E., and Leung, C. K. Y. (2017), {A coupled diffusion-mechanical model with boundary element method to predict concrete cover cracking due to steel corrosion}, \textit{Corrosion Science}, \textbf{126}, {180--196}.

\bibitem[\protect\citeauthoryear{Coelho \textit{et~al.}}{2017}]{coelho2017mechanical}
Coelho, K. O. and Leonel, E. D. and Fl{\'o}rez-L{\'o}pez, J. (2017), {The Mechanical Behavior Modeling of Reinforced Concrete Structures by the Lumped Damage Model.} In: Robinson, S., \textit{Reinforced Concrete: Design, Performance and Applications}, Nova Science Publishers, Inc, Hauppauge, New York, USA.

\bibitem[\protect\citeauthoryear{Estes and Frangopol}{1998}]{estes1998relsys}
Estes, A. C and Frangopol, D. M. (1998), {RELSYS: A computer program for structural system reliability},  \textit{Structural engineering and mechanics}, \textbf{8}, 901--919.

\bibitem[\protect\citeauthoryear{Fl\'orez-L\'opez}{1998}]{florez1998frame}
Fl\'orez-L\'opez, J. (1998), {Frame analysis and continuum damage mechanics},  \textit{European Journal of Mechanics-A/Solids}, \textbf{17}, 269--283.

\bibitem[\protect\citeauthoryear{Goctermann}{2000}]{goctermann2000dura}
Goctermann, P. (2000), {DuraCrete probabilistic performance based durability design of concrete structure: general guidelines for durability design and redesign},  \textit{Report No. BE9521347/R14}.

\bibitem[\protect\citeauthoryear{Griffith}{1921}]{griffith1921vi}
Griffith, A. A. (1921), {VI. The phenomena of rupture and flow in solids}, \textit{Philosophical transactions of the royal society of london. Series A, containing papers of a mathematical or physical character}, \textbf{221}, 163--198.

\bibitem[\protect\citeauthoryear{Huang \textit{et al.}}{2020}]{huang2020stochastic}
Huang, L., Jin, X., Fu, C., Ye, H., and Dong, X. (2020), {Stochastic characteristics of reinforcement corrosion in concrete beams under sustained loads}, \textit{Computers and Concrete}, \textbf{25}(5), {447--460}.

\bibitem[\protect\citeauthoryear{Jia \textit{et al.}}{2020}]{jia2020stochastic}
Jia, G., Gardoni, P., Trejo, D, and Mazarei, V. (2020), {Stochastic Modeling of Deterioration and Time-Variant Performance of Reinforced Concrete Structures under Joint Effects of Earthquakes, Corrosion, and ASR}, \textit{Journal of Structural Engineering}, \textbf{147}(2), {04020314}.

\bibitem[\protect\citeauthoryear{Kiani and Shodja}{2009}]{kiani2012response}
Kiani, K. and Shodja, H. (2012), {Response of reinforced concrete structures to macrocell corrosion of reinforcements. Part II: After propagation of microcracks via a numerical approach},  \textit{ Nuclear Engineering And Design}, \textbf{242}, 7--18.

\bibitem[\protect\citeauthoryear{Liberati \textit{et~al.}}{2014}]{liberati2014nonlinear}
Liberati, E. A. P, Nogueira C. G., and Leonel, E. D. (2014), {Nonlinear formulation based on FEM, Mazars damage criterion and Fick’s law applied to failure assessment of reinforced concrete structures subjected to chloride ingress and reinforcements corrosion},  \textit{Engineering Failure Analysis}, \textbf{46}, 247--268.

\bibitem[\protect\citeauthoryear{Marante and Fl\'orez-L\'opez}{2003}]{marante2003three}
Marante, M. E. and Fl\'orez-L\'opez, J. (2003), {Three-dimensional analysis of reinforced concrete frames based on lumped damage mechanics},  \textit{International Journal of Solids and Structures}, \textbf{40}, 5109--5123.

\bibitem[\protect\citeauthoryear{Oudah and Norlander}{2019}]{oudah2017design}
Oudah, F. and Norlander, F. (2017), {Design philosophy and requirements of granular wear surface thickness for bridges subjected to extreme truck load}, \textit{Canadian Journal of Civil Engineering}, \textbf{44}, 727--735.

\bibitem[\protect\citeauthoryear{Otieno \textit{et~al.}}{2012}]{otieno2012prediction}
Otieno, M., Beushausen, H. and Er, M. (2012), {Prediction of corrosion rate in reinforced concrete structures--a critical review and preliminary results},  \textit{Materials And Corrosion}, \textbf{63}, 777--790.

\bibitem[\protect\citeauthoryear{Papadakis \textit{et al.}}{1992}]{papadakis1992effect}
Papadakis, V.G., Fardis, M.N. and Vayenas, C.G. (1992), {Effect of composition, environmental factors and cement-lime mortar coating on concrete carbonation}, \textit{Materials and structures}, \textbf{25}(5), {293--304}.

\bibitem[\protect\citeauthoryear{Papadakis \textit{et al.}}{1996}]{papadakis1996mathematical}
Papadakis, V.G., Roumeliotis, A.P., Fardis, M.N. and Vagenas, C.G. (1996), {Mathematical modelling of chloride effect on concrete durability and protection measures}, In: Dhir, R. K. and Jones, M. R., \textit{Concrete repair, rehabilitation and protection}, E \& FN Spon, London, UK.

\bibitem[\protect\citeauthoryear{Pellizzer \textit{et~al.}}{2015}]{pellizzer2015mechanical}
Pellizzer, G. P., Leonel, E. D. and Nogueira, C. G. (2015), {Infuence of reinforcement's corrosion into hyperstatic reinforced concrete beams: a probabilistic failure scenarios analysis},  \textit{RIEM-IBRACON Structures and Materials Journal}, \textbf{8}, 479--490.

\bibitem[\protect\citeauthoryear{Pellizzer \textit{et al.}}{2020}]{pellizzer2020time}
Pellizzer, G. P., Kroetz, H. M., Leonel, E. D. and Beck, A. T. (2020), {Time-dependent reliability of reinforced concrete considering chloride penetration via boundary element method}, \textit{Latin American Journal of Solids and Structures}, \textbf{17}(8),{1--17}.

\bibitem[\protect\citeauthoryear{Santoro and Kunnath}{2013}]{santoro2013damage}
Santoro, M. G. and Kunnath, S. K. (2013), {Damage-based RC beam element for nonlinear structural analysis}, \textit{Engineering structures}, \textbf{49},{733--742}.

\bibitem[\protect\citeauthoryear{Shafikani and Chidiac}{2019}]{shafikhani2019quantification}
Shafikhani, M. and Chidiac, S.E. (2019), {Quantification of concrete chloride diffusion coefficient--A critical review}, \textit{Cement and Concrete Composites}, \textbf{99}, 225--250.

\bibitem[\protect\citeauthoryear{Shah and Bishnoi}{2018}]{shah2018carbonation}
Shah, V. and Bishnoi, S. (2018), {Carbonation resistance of cements containing supplementary cementitious materials and its relation to various parameters of concrete}, \textit{Construction and Building Materials}, \textbf{178}, {219--232}.

\bibitem[\protect\citeauthoryear{Stewart}{2004}]{stewart2004spatial}
Stewart, M. G. (2014), {Spatial variability of pitting corrosion and its influence on structural fragility and reliability of RC beams in flexure},  \textit{Structural Safety}, \textbf{26}, 453--470.

\bibitem[\protect\citeauthoryear{Val \textit{et~al.}}{2000}]{val1998effect}
Val, D. V., Stewart, M. G. and Melchers, R. E. (1998), {Effect of reinforcement corrosion on reliability of highway bridges}, \textit{Engineering structures}, \textbf{20}, 1010--1019.

\bibitem[\protect\citeauthoryear{Vecchio and Collins \textit{et~al.}}{1986}]{vecchio1986modified}
Vecchio, F. J. and Collins, M. P. (1986), {The modified compression-field theory for reinforced concrete elements subjected to shear}, \textit{ACI Journal}, \textbf{83}, 219--231.

\bibitem[\protect\citeauthoryear{Vennesland \textit{et~al.}}{2009}]{vennesland2013recommendation}
Vennesland, O., Climent, M.-{\'A}. and Andrade, C. (2013), {Recommendation of RILEM TC 178-TMC: Testing and modelling chloride penetration in concrete}, \textit{Materials and Structures}, \textbf{46}, 337--344.

\bibitem[\protect\citeauthoryear{Vu and Stewart}{2000}]{vu2000structural}
Vu, K. A. T. and Stewart, M. G. (2000), {Structural reliability of concrete bridges including improved chloride-induced corrosion models}, \textit{Structural safety}, \textbf{22}, 313--333.

\bibitem[\protect\citeauthoryear{Tuutti}{1982}]{tuutti1982corrosion}
Tuutti, K. (1982), \textit{Corrosion of steel in concrete},  Cement-och betonginst, Stockholm, Sweden.

\bibitem[\protect\citeauthoryear{Ueda and Takewaka}{2007}]{ueda2007performance}  
Ueda, T. and Takewaka, K. (2007), {Performance-based standard specifications for maintenance and repair of concrete structures in Japan}, \textit{Structural Engineering International}, \textbf{17}, 359--366.

\bibitem[\protect\citeauthoryear{Xiao \textit{et al.}}{2012}]{xiao2012fem}
Xiao, J., Ying, J. and Shen, L. (2012), {FEM simulation of chloride diffusion in modeled recycled aggregate concrete}, \textit{Construction and Building Materials}, \textbf{29}, {12--23}.

\bibitem[\protect\citeauthoryear{Xu \textit{et al.}}{2018}]{xu2018numerical}
Xu, F., Xiao, Y., Wang, S., Li, W., Liu, W. and Du, D. (2018), {Numerical model for corrosion rate of steel reinforcement in cracked reinforced concrete structure}, \textit{Construction and Building Materials}, \textbf{180}, {55--67}.

\bibitem[\protect\citeauthoryear{Yu \textit{et al.}}{2021}]{yu2021integrated}
Yu, Y., Gao, W., Castel, A., Chen, X., and Liu, A. (2021), {An integrated framework for modelling time-dependent corrosion propagation in offshore concrete structures}, \textit{Engineering Structures}, \textbf{228}, {111482}.
\end{thebibliography}
\end{document}